\documentclass[11pt]{article}

\usepackage{color}
\usepackage{fullpage}

\usepackage{amssymb}
\usepackage{amsbsy}
\usepackage{amsthm}
\usepackage{enumerate}
\usepackage{framed}
\usepackage[utf8]{inputenc}
\usepackage[T1]{fontenc}
\usepackage{imakeidx}

\newcommand{\ad}{{\rm ad}}

\newcommand{\real}{{\mathbb R}}
\newcommand{\complex}{{\mathbb C}}

\newcommand{\wh}[1]{\widehat{#1}}

\newcommand{\fract}[2]{\textstyle{\frac{#1}{#2}}}
\newcommand{\ip}[2]{\langle{#1}|{#2}\rangle}

\newcommand{\LieE}{{{\mathfrak e}}}
\newcommand{\LieF}{{{\mathfrak f}}}
\newcommand{\LieG}{{{\mathfrak g}}}

\newcommand{\LieT}{{{\mathfrak t}}}
\newcommand{\so}{{{\mathfrak {so}}}}
\newcommand{\spin}{{{\mathfrak {spin}}}}
\newcommand{\symp}{{{\mathfrak {sp}}}}
\newcommand{\su}{{{\mathfrak {su}}}}

\newcommand{\schur}{{\mathbb {S}}}

\newcommand\rmap[1]{\stackrel{#1}\longrightarrow}

\begin{document}

\begin{center}
{\LARGE{Tensor powers of adjoint representations of classical Lie groups}}

\bigskip
{K.C.\ Hannabuss

\medskip
Mathematical Institute, ROQ, Woodstock Road, 
Oxford OX2 6GG, UK\\ 
\smallskip
{\it Email: hannabus@maths.ox.ac.uk}}

\end{center}

\begin{abstract}
Exploiting particular features of classical groups, simple constructions are given for the irreducible constituents of the tensor square of the adjoint modules  and the leading terms in higher tensor powers. This provides an independent confirmation of Vogel's general formulae, and alternative approach to that in some recent papers.
\end{abstract}


\section*{Introduction}

To lighten notation we use the terminology of $G$-modules 
rather than representations of the group $G$. (A $G$-module is a vector space $M$ equipped with a map $G\times M \to M$ ($(g,m) \mapsto g.m$) which satisfies linearity in $M$ and $g.(h.m) = (g.h).m$, and corresponds to a representation $g: m \mapsto g.m$ on $M$, and, for consistency,  the adjoint module of a Lie group $G$ is the $G$-module corresponding to the adjoint representation.) Functors such as symmetric or exterior tensor powers provide new $G$-modules.
The $k$-th Cartan power $M^{(k)}$ of an irreducible module $M$ is the irreducible with highest weight $k$ times that of $M$. The dual of $M$ will be denoted by $\wh{M}$.

We also write $\LieG$ for the complex Lie algebra of $G$, and use it for both the algebra and the adjoint $G$-module itself.
The adjoint module of any classical group is easily linked to some tensor product of its natural module with itself or its dual.
This quickly leads to simple constructions of the irreducible constituents of tensor powers of their 
adjoint modules which may provide a useful complement to the direct approach to that problem in \cite{IK1,IK2}.
(Our constructions are simple enough that they may be well-known, but the author could not locate any reference.)
After  describing the irreducible decomposition of the tensor square of the adjoint module, we give the most interesting components of its higher tensor powers, whose dimensions match those calculated in \cite{LM}.
One motivation for this came from the deep results of Vogel \cite{V,V2}, which simultaneously classify all simple Lie algebras and decompose the tensor squares of their adjoint modules in a uniform way, exploiting an extension of diagrammatic methods developed for finding invariants in knot theory. The present paper does not address the classification, but does describe some straightforward explicit decompositions of the tensor squares and some higher powers. Another motivation came from the 
alternative definition of  highest weight vectors, \cite{KCH82,WL82}, a new proof of which is given in  
Appendix B.

The tensor product $\LieG\otimes\LieG$ decomposes into the direct sum of the symmetric and antisymmetric (exterior) tensor products  $\schur^2\LieG = \LieG\otimes_S\LieG$ and $\wedge^2\LieG = \LieG\wedge\LieG$, respectively.
The symmetric tensor product contains the Cartan square $\LieG^{(2)}$, and it clearly also contains a trivial submodule generated by the polarised Casimir operator.
For some small Lie algebras these two suffice, and 
we have $\schur^2\LieG \cong \LieG^{(2)} \oplus \complex$, but these two submodules are usually insufficient to span the whole space.
The Lie bracket provides a projection $X\wedge Y \mapsto [X,Y]$ from $\wedge^2\LieG$ onto $\LieG$, so that $\LieG$ is always a submodule of $\wedge^2\LieG$, but usually some more submodules also appear.


\section{The tensor square of the adjoint module of a unitary group}

The special unitary groups ${\rm SU}(n)$ with $n\geq 3$ have a cubic as well as a quadratic Casimir operator, which, together, provide maps from $\schur^2 \LieG \to \wh{\LieG} \cong \LieG$, and its image provides a copy of $\LieG$ as a submodule of $\schur^2\LieG$.
(For example, when $G = {\rm SU}(3)$, $\schur^2 \LieG$ splits into the direct sum of the 27-dimensional Cartan square $\LieG^{(2)}$,  a copy of the 8-dimensional $\LieG$, and the trivial module $\complex$.)

The group $G = {\rm SU}(n)$ has its natural module  $V= \complex^n$, as do ${\rm U}(n)$ and the complexifications ${\rm SL}(n,\complex)$ and ${\rm GL}(n,\complex)$, and the linear transformations $\hom(V)$ on $V$ can be identified with the Lie algebra of the last group. Under the conjugation action of ${\rm SU}(n)$ and ${\rm U}(n)$  $\hom(V)$  decomposes into the direct sum of the trace zero matrices, which can be identified with the complex Lie algebra $\LieG = \su(n)$  and  the multiples of the identity operator: 
$\hom(V) \cong \LieG \oplus 1$.
The symmetric and antisymmetric tensor products are therefore
\begin{eqnarray*} 
\schur^2\hom(V) &\cong& \schur^2 \LieG \oplus \LieG \oplus \complex1 \cong \schur^2 \LieG \oplus \hom(V)\\
\wedge^2 \hom(V) &\cong& \wedge^2\LieG \oplus \LieG,
\end{eqnarray*}
where the copy of $\hom(V)$ in the symmetric product is  
$\Delta(X) = \{(X\otimes 1+ 1\otimes X): X\in \hom(V)\}$,
and of $\LieG$ in  the antisymmetric product is $\{(X\otimes 1- 1\otimes X): X\in \LieG\}$.

We can also write $\hom(V) \cong \wh{V}\otimes V$ (where $\wh{V}$ is the dual of $V$, and then
\begin{eqnarray*}
\schur^2\hom(V) 
&\cong& (\schur^2 \wh{V} \otimes \schur^2V) \oplus (\wedge^2\wh{V} \otimes \wedge^2V) 
\cong \hom(\schur^2 V) \oplus \hom(\wedge^2 V)\\
\wedge^2 \hom(V) 
&\cong& (\schur^2 \wh{V}) \otimes (\wedge^2 V) \oplus (\wedge^2 \wh{V}) \otimes (\schur^2V) 
\cong \hom(\schur^2V,\wedge^2V) \oplus \hom(\wedge^2V, \schur^2V). 
\end{eqnarray*}
The two summands in the latter case are automatically dual to each other.

These simple observations already provide considerable information.
When we choose an element of a maximal abelian subgroup $T$ of $\hom(V) \sim {\mathfrak {gl}}(V)$, whose action on $V$ has eigenvalues $\alpha_1, \dots, \alpha_n$ then the eigenvalues of the action on $\schur^2V$ are $\{\alpha_j\alpha_k: j\leq k\}$, and on $\wedge^2V$ are $\{\alpha_j\alpha_k: j<k\}$ . 
The highest weight on $\hom(\schur^2V)$ is $\alpha_1^2\alpha_n^{-2}$ or, identifying the dual Lie algebra $\wh{\LieT}$ with the subspace of $\real^n$ whose components sum to 0, the highest weight is identified with  $(2,0,\ldots,0,-2) \in \wh{\LieT}$, which we usually abbreviate to  $(2,0^{n-2},-2)$. This is the highest weight of the Cartan square $\LieG^{(2)}$, since $\LieG$ has highest root is $(1,0,\ldots,0,-1)$.
The dimension of the irreducible module with this highest weight is 
$$
\dim[\LieG^{(2)}] = \fract14n^2(n-1)(n+3) = [\fract12n(n+1)]^2 - n^2.
$$ 
Since $\dim[\hom(\schur^2V)] = (\frac12n(n+1))^2$, and $\dim[\hom(V)] = n^2$,  we see
 that there is no room for more irreducibles or multiplicities so $\hom(\schur^2V) = \LieG^{(2)}\oplus \hom(V)$.
 Throughout our discussions the dimensions can be found using Weyl's dimension formula (see Appendix A),  or Stanley's hook content formula, \cite[7.21.2]{S}, \cite[4.6]{FH}.

On $\hom(\wedge^2V)$  (provided that $n \geq 4$) the highest weight is $(1,1,0,\ldots,0,-1,-1)$,  which we abbreviate to  $(1^2,0^{n-4},(-1)^2)$,  
and we shall denote the module with this highest weight by $\LieG^{(1^2)}$. 
Its dimension is given by
$$
\dim[\LieG^{(1^2)}] = \fract14n^2(n+1)(n-3) = [\fract12n(n-1)]^2 - n^2,
$$ 
from which we deduce the decomposition
$\hom(\wedge^2V) = \LieG^{(1^2)}\oplus \hom(V)$.

In the case of $\wedge^2 \hom(V)$, we noted that the direct summands 
$\hom(\schur^2V,\wedge^2V)$ and $\hom(\wedge^2V, \schur^2V)$ are dual and each has dimension
$$
\dim[\hom(\wedge^2V,\schur^2V)] = \left[\fract12n(n+1)\right]\left[\fract12n(n-1)\right]
= \fract14[n^2(n^2-1)].
$$
The highest weight in $\hom(\wedge^2V, \schur^2V)$ is $\alpha_1^2(\alpha_{n-1}\alpha_n)^{-1}$ and in 
$\hom(\schur^2V,\wedge^2V)$ is $\alpha_1\alpha_2\alpha_n^{-2}$. Both the irreducibles with these highest weights have dimension
$\frac14(n^2-1)(n^2-4)$, and we label them as $\LieG^{(1^2,2)}$, $\LieG^{(2,1^2)}$, respectively.
In this case the homomorphism spaces contain only a copy of $\LieG$, dimension $n^2-1$ which is just right to be the its orthogonal complement of theses two irreducibles.

Already  $\LieG^{(2)}$,  $\LieG^{(1^2)}$, $\LieG^{(1^2,2)}$, and $\LieG^{(2,1^2)}$ have exactly the dimensions of the irreducible submodules as calculated by Vogel and others \cite{V2,LM,IK1,IK2}.  
Pulling everything together, we arrive at
\begin{eqnarray*}
\schur^2 \LieG &\cong& \LieG^{(2)} \oplus \LieG^{(1^2)} \oplus \LieG\oplus \complex1,\\
\wedge^2 \LieG &\cong& \LieG^{(1^2,2)} \oplus \LieG^{(2,1^2)} \oplus \LieG.
\end{eqnarray*}
However, we should note that for $\dim(V)\leq 4$ some of these submodules are absent: 
for example, $\wedge^2V$ is the dual of $V$ when $\dim(V) =3$, and so $\dim(\LieG^{(1^2)}) = \dim(\hom(\wedge^2V) - \dim(\hom(V)) = 0$.

\medskip
\noindent{\it Remark.}
Surprisingly, the duality of components in $\wedge^2V$ can be extended to $\schur^2V$ by using 
virtual modules.
It is easy to check the following identity for modules $V$ and $W$:
$$
\wedge^2(V) \oplus \schur^2(V\oplus W) 
\cong V\otimes (V \oplus W) \oplus\schur^2(W),
$$
and on setting $W$ equal to the virtual module $-V$ this reduces to 
$\wedge^2(V) \cong \schur^2(-V)$. On the other hand $\hom(-V) \cong \hom(V)$, so that the connections with $\LieG$ are preserved, but, when $V \mapsto -V$, $\hom(\schur^2V)$ and $\hom(\wedge^2V)$ are interchanged, as are  the dual pair $\hom(\wedge^2V,\schur^2V)$ and $\hom(\schur^2V,\wedge^2V)$.
Now $\hom(\schur^2V)$ and $\hom(\wedge^2V)$ each contain a copy of $\hom(V)$ complementing 
$\LieG^{(2)}$ and $\LieG^{(1^2)}$, respectively, so that changing $V$ to $-V$ interchanges $\LieG^{(2)}$ and 
$\LieG^{(1^2)}$.
\medskip


\section{Higher tensor powers of unitary group adjoint modules}

It is easy to extend some of the analysis to higher tensor powers $\schur^k\LieG$, but important changes occur for $k>2$. 
As before, $\schur^k\hom(V)$ still contains $\hom(\schur^kV)$ and $\hom(\wedge^kV)$, and
$\wedge^k\hom(V)$ still contains $\hom(\wedge^kV,\schur^kV)$ and $\hom(\schur^kV,\wedge^kV)$, 
and we now have
$$
\schur^k\kern-2pt\hom(V)  \cong \schur^k\LieG \oplus \schur^{k-1}\kern-2pt\hom(V) \qquad\wedge^k\hom(V) \cong \wedge^k\LieG \oplus \wedge^{k-1}\LieG.
$$  
However, symmetric and antisymmetric tensors are no longer sufficient to to handle higher tensor powers $\bigotimes^k\!\hom(V)$, since for any partition $\mu$ of $k$ a Schur functor $\schur_\mu$ provides a submodule
$\schur_\mu V$ of $\bigotimes^kV$; these include $\schur_{(k)} = \schur^k$, and $\schur_{(1^k)} = \wedge^k$. In contrast to the tensor square where  $\schur_{(2)} = \schur^2$, and $\schur_{(1^2)} = \wedge^2$ sufficed,  the additional partition $3 = 2+1$ means that $\schur_{(2,1)}$ is needed in the tensor cube. 
For $k > n$ many partitions $\mu$ of $k$ with columns longer than $n$ will give $\schur_\mu V = 0$.

Schur--Weyl duality leads us to expect the decompositions to involve terms like $\hom(\schur_\nu V,\schur_\mu V)$.  
The action of the symmetric group $S_k$ on $\schur^k\hom(V) = \schur_{(k)}\hom(V)$ is trivial and the decomposition requires only terms where $\mu = \nu$, giving, \cite[Ex. 6.11(a)]{FH}, whilst exterior powers $\wedge^k \sim \schur_{(1^k)}$ corresponding to the alternating module require $\mu$ and $\nu = \mu^\prime$ to be dual partitions. (In $\mu^\prime$  the rows and columns of $\mu$ are interchanged.) In other words, we have 
$$
\schur^k\hom(V) \cong \oplus_{\mu}  \hom(\schur_\mu V), \qquad
\wedge^k\hom(V) \cong \oplus_{\mu}  \hom(\schur_{\mu^\prime}V,\schur_{\mu}V),
$$
with dimensions, expressed in terms of $D_\mu = \dim[\schur_\mu(V)]$, as 
$$
\dim[\schur^k\hom(V)] = \sum_{\mu}  D_\mu^2, \qquad
\dim[\wedge^k\hom(V)]  = \sum_{\mu} D_{\mu^\prime}D_\mu.
$$

Working with virtual modules we find that
$\schur^k\hom(V) -\wedge^k\hom(V)$ is isomorphic to 
\begin{eqnarray*}
&&\oplus_{\mu}  [\hom(\schur_\mu V)
- \hom(\schur_{\mu^\prime}V,\schur_{\mu}V)]\\
&&\cong \oplus_{\mu}^\prime[\hom(\schur_\mu V) +\hom(\schur_{\mu^\prime}V)
- \hom(\schur_{\mu^\prime} V,\schur_{\mu}V)- \hom(\schur_{\mu}V,\schur_{\mu^\prime}V)]\\
&&\cong  \oplus_{\mu}^\prime[\hom(\schur_\mu V-\schur_{\mu^\prime}V)],
\end{eqnarray*}
where $\oplus_{\mu}^\prime$ includes only one of the two terms involving $\mu$ and $\mu^\prime$.
(All self-dual partitions cancel in this formula (so that, for example, 
$\schur^3\hom(V) -\wedge^3\hom(V) \cong  \hom(\schur^3V -\wedge^3V)$), and only half the others need be considered since $\mu$  and $\mu^\prime$ give equivalent contributions, 
providing a useful way of deducing results about $\schur^k\hom(V)$ from those for 
$\wedge^k\hom(V)$.)  
The dimensions are related by
\begin{eqnarray*}
\dim[\schur^k\hom(V) -\wedge^k\hom(V)]
&=&\sum_{\mu}{}^\prime  [D_\mu^2 - D_{\mu^\prime} D_\mu]\\
&=& \sum_{\mu}  \fract12[D_\mu^2 + D_{\mu^\prime}^2 - 2D_{\mu^\prime} D_\mu]\\
&=&  \sum_{\mu} \fract12[D_\mu - D_{\mu^\prime}]^2,
\end{eqnarray*}
and when $k=3$ we see that 
$\dim[\schur^3\hom(V) -\wedge^3\hom(V)] =  \dim[\schur^3V -\wedge^3V]^2$.

There is a direct way to embed $\hom(\schur^{k-1}V)$ into $\hom(\schur^kV)$ which generalises the case of $k=2$. We first send $A \in \hom(\schur^{k-1}V)$ to $A\otimes I$ and then symmetrise this by averaging over the cyclic subgroup generate by the cyclic permutation $\tau =(12\ldots k)$ which projects $A\otimes I$ to $\sum_{j=1}^{k} \tau^j(A\otimes I)\tau^{-j}$.
Since $A$ is symmetrised this is symmetric under $S_k$. 
Taking care to insert factors of $-1$ for odd permutations, there is an analogous map 
$\hom(\wedge^{k-1}V) \to \hom(\wedge^kV)$.

To keep things straightforward and avoid the other partitions, we shall only consider the modules 
$\hom(\schur^kV), \hom(\wedge^kV) \subseteq \schur^k\hom(V)$,
which suffice to give us information comparable with \cite{LM}. 
(Moreover, for special unitary groups the isomorphism $\wedge^kV \cong \wedge^{n-k}V$ means that it is sufficient to consider only $k \leq \frac12 n$ in this case.)

The complement of $\hom(\schur^{k-1}V)$ in $\hom(\schur^kV)$ has dimension
\begin{eqnarray*}
{{n+k-1}\choose{k}}^2 - {{n+k-2}\choose{k-1}}^2 &=&  \frac{(n+2k-1)}{(n-1)}{{n+k-2}\choose{k}}^2.
\end{eqnarray*}
Similarly, the complement of $\hom(\wedge^{k-1}V)$ in $\hom(\wedge^kV)$ has dimension
\begin{eqnarray*}
{{n}\choose{k}}^2 - {{n}\choose{k-1}}^2  &=&  \frac{(n-2k+1)}{(n+1)}{{n+1}\choose{k}}^2.
\end{eqnarray*}
Our observation that we could use $k\leq \frac12 n$ removes the obvious problems with this formula when $k \leq n < 2k$.

It is shown in Appendix A that both of these are irreducible. The former is the Cartan power 
$\LieG^{(k)}$, and the latter will be denoted by $\LieG^{(1^k)}$.
Using the Vogel parameters for special linear groups, ($\alpha = -2 = -\beta$, $\gamma = n = t$), their dimensions agree with those in \cite[Thm 1.2]{LM} : $\dim[\LieG^{(k)}]= Y_k(\alpha)$ 
and $\dim[\LieG^{(1^k)}] = Y_k(\beta)$.

In the study of Casimir operators in the next section, we shall need the highest weights of the various modules and these are easily obtained.
The highest weight on $\hom(\schur^kV)$ is $(k,0^{n-2},-k)$, which is the $k$-th power of the highest root, and so is the highest weight of the Cartan power $\LieG^{(k)}$.

 
\section{The Casimir operators}
 
We denote the nondegenerate symmetric bilinear Killing form on $\LieG$ by 
$\kappa: \LieG^{\otimes 2} \to \complex$. The non-degeneracy permits us to  form its dual 
$\wh{\kappa}: \complex \to \schur^2\LieG$, satisfying $\kappa^{\otimes2}(\wh{\kappa},X\otimes Y) = \kappa(X,Y)$ for all $X, Y \in \LieG$, and we shall identify $\wh{\kappa}$ with the image of 1 
in $\schur^2\LieG$ and interpret it as the polarised Casimir element.
In condensed Sweedler notation (where summation is suppressed) we write $\wh{\kappa} = \wh{\kappa}_1\otimes \wh{\kappa}_2$.  
Then, by definition, we have, for all $X, Y \in \LieG$,
$\kappa(\wh{\kappa}_1,X)\kappa(\wh{\kappa}_2,Y) = \kappa(X,Y)$, giving the expansion formula 
$\kappa(\wh{\kappa}_1,X)\wh{\kappa}_2 = X$, and, since $\wh{\kappa}$ is symmetric, a similar expansion with $\wh{\kappa}_1$ and $\wh{\kappa}_2$ interchanged. We can also express this as a formula for the identity operator $I = \wh{\kappa}_1 \kappa(\wh{\kappa}_2,\cdot) = \wh{\kappa}_2 \kappa(\wh{\kappa}_1,\cdot)$.

(We can normalise the Killing form to be consistent with earlier conventions, so that $\wh{\kappa}$ is the normal Euclidean inner product on the weight vectors $\lambda$, and all roots such as $\lambda =(1,0,0\ldots,0,-1)$ give $\wh{\kappa}(\lambda,\lambda) = 2$.  This avoids unwanted factors, which would otherwise complicate formulae, though this normalisation is not universally used (in  particular, \cite{IK1,IK2} use very different conventions).

The (quadratic) Casimir element $C$ in the enveloping algebra of $\LieG$ is obtained by applying the algebra multiplication $m$ to get $C = m\circ\wh{\kappa}$. 
We can also form a quadratic Casimir element $C^\Delta$ for tensor products by replacing each $X \in \LieG$ appearing in $C$ by its comultiplication $\Delta X = X\otimes 1 + 1\otimes X$ and, we may then check that 
$$
C^\Delta = C\otimes 1 + 1\otimes C + 2\wh{\kappa},
$$
which displays $\wh{\kappa}$ as a polarisation of the Casimir operator (It is sometimes referred to as {\lq split\rq\ rather than polarised, \cite{IK1,IK2}.)  When $V$ is an irreducible module of $G$,  Schur's Lemma gives $C = c1$ for some $c \in \complex$.  (We generally arrange things so that $C$ is self-adjoint and then $c \in \real$).  In that situation $C^\Delta$ and $\wh{\kappa}$ contain the same information and the eigenvalues $\{c_j\}$ of $C^\Delta$, and $\{k_j\}$ of $\wh{\kappa}$ are related by $c_j = 2(c + k_j)$. 

The quadratic Casimir acts on any $G$-module $M$, whilst $C^\Delta$ and the polarised Casimir $\wh{\kappa} \in \LieG \otimes \LieG$ automatically act on the tensor product  $M_1\otimes M_2$ of two $G$-modules,  $M_1$ and $M_2$.
In contrast to the Casimir operator, $C^\Delta$ whose eigenvalues on irreducible submodules can be found from the highest weights by the usual formula,  we see that $\wh{\kappa}$ is sensitive to the context of a submodule. 

When $M_1 = M_2= M$ we have the identity
$\wh{\kappa}= \fract12(C^\Delta - C\otimes 1 - 1\otimes C)$ on $M\otimes M$, and,  for a module  with highest weight $\lambda$ the polarised Casimir acting on the submodule $M^{(2)}$ of $M\otimes M$ with highest weight $2\lambda$ is
$$
\fract12\bigl[\wh{\kappa}(2\lambda,2\lambda+2\delta) - 2\wh{\kappa}(\lambda,\lambda+2\delta)\bigr]
= \wh{\kappa}(\lambda,\lambda).
$$

The number $N$ of distinct eigenvalues $\{c_j\}$ of $C^\Delta$ is finite, and so the operator satisfies the minimal equation 
$$
\prod_{j = 1}^N(C^\Delta - c_j1) = 0.
$$
The projection onto the $c_r$-eigenspace is given by 
$$
P_r = \prod_{j \neq r}\frac{(C^\Delta - c_j1)}{(c_r - c_j1)}.
$$
The eigenspace is the image of $P_r$ which is also $\ker(C^\Delta - c_j1)$, and its orthogonal complement is $(C^\Delta - c_j1)(V\otimes V)$. There are clearly analogous formulae in terms of $\wh{\kappa}$.

The sum of the positive roots takes the form $2\delta = (n-1,n-3,n-5,\ldots, -(n-1))$ for the unitary groups.
With our rescaling of $\wh{\kappa}$ and identification of the dual of $\LieT$ it reduces to the Euclidean inner product
$\wh{\kappa}(\lambda,\mu) =  \sum_j \lambda_j\mu_j$, 
\cite[Eq. (15.4)]{FH}.
For the adjoint module we find that $\wh{\kappa}(\lambda,\lambda) = 2$ and $\wh{\kappa}(\lambda,2\delta) = 2(n-1)$ so that the eigenvalue of the quadratic Casimir is $2(n-1) + 2 = 2n$.
For the Cartan powers $\LieG^{(k)}$ we get
$$
\wh{\kappa}(k\lambda,k\lambda + 2\delta) = 2k^2 + 2k(n-1) = 2k(n+k-1).
$$

Finally, for $\LieG^{(1^k)}$ we have 
$$
\wh{\kappa}(\lambda,\lambda+2\delta) = 2k + 2\sum_{j=1}^k (n+1-2j) = 2k(n+2) -2k(k+1) 
= 2k(n-k+1).
$$

On the submodules $\LieG^{(1^k,k)}$ and $\LieG^{(k,1^k)}$, we find that the Casimir operator takes the value $2kn$.

To obtain the eigenvalues of the polarised Casimir operators on $\schur^2V$ we simply take half the quadratic Casimir value when $k=2$, and subtract its value $2n$ in the adjoint module to get
the eigenvalues of the polarised Casimir operator $\wh{\kappa}_\schur$ on each non-trivial submodule.
\begin{eqnarray*}
{\LieG^{(2)}}:  2(n+1) - 2n = 2,&\qquad&  
{\LieG^{(1^2)}}:  \frac12(4(n-1)) - 2n = -2, \\
\LieG: \frac12(2n) - 2n = -n, &\qquad&
\LieG^{(2,1^2)} =
\LieG^{(1^2,2)}:   2n - 2n = 0.
\end{eqnarray*}
The eigenvalues for the nontrivial submodules $\LieG^{(2)}$, $\LieG^{(1^2)}$, and $\LieG$ in the symmetric part are $2$, $-2$, and $-n$, which are just the additive inverses of the Vogel parameters $\alpha$, $\beta$, and $\gamma$ for the linear groups, \cite{V,V2,LM,IK1,IK2}.
The  trivial submodule of the adjoint module is in the kernel of the quadratic Casimir operator, so  $\lambda =0$, and the polarised Casimir takes the value $-2n$. 

The eigenvalues of the quadratic Casimir operators are so easily calculated that it is easier to find the eigenvalues of the polarised Casimir this way, but direct calculations are possible. They can be simplified using the fact that $\schur^2\LieG$ is spanned by elements $X^{\otimes2}$ for $X\in \LieG$, and the action of $\wh{\kappa}$ is given by
$\wh{\kappa}\cdot X^{\otimes2} = [\wh{\kappa}_1,X]\otimes_S [\wh{\kappa}_2,X]$,
and $\kappa^{\otimes2}(Y^{\otimes2},\wh{\kappa}\cdot\! X^{\otimes2})=\kappa([X,Y],[X,Y])$.

Since we are in the semi-simple situation and the only eigenvalues of the polarised Casimir operator $\wh{\kappa}_\schur$ on $\schur^2\LieG$ are $\pm 2$, $-n$ and $-2n$ we deduce that it satisfies the minimal equation (usually known in this context as the characteristic equation)
$$
(\wh{\kappa}_\schur^2 -4)(\wh{\kappa}_\schur+n)(\wh{\kappa}_\schur+2n) 
= 0.
$$
The Lagrange interpolation formula gives the projections $P_{-2}$, $P_2$, $P_{-n}$ and $P_{-2n}$ onto the different irreducible summands, such as
$$
P_ {\pm2}= \frac{(\wh{\kappa}_\schur\pm2)(\wh{\kappa}_\schur+n)(\wh{\kappa}_\schur+2n)}{\pm4(n\pm2)(2n\pm2)}, \quad
P_ {-n}= \frac{(\wh{\kappa}_\schur^2 -4)(\wh{\kappa}_\schur+2n)}{n(n^2 - 4)},\quad
P_ {-2n}= -\frac{(\wh{\kappa}_\schur^2 -4)(\wh{\kappa}_\schur+n)}{4n(n^2 - 1)},
$$
and $\wh{\kappa}_\schur =  2 P_{+2} -2 P_{-2} - n P_{-n} -2n P_{-2n}$.

On $\wedge^2\LieG$ the polarised Casimir $\wh{\kappa}_\wedge$ has only the eigenvalues $-n$ and $0$, giving 
$$
\wh{\kappa}_\wedge(\wh{\kappa}_\wedge+n)=  0, 
$$
and projections $P_{-n} = -n^{-1}\wh{\kappa}_\wedge$, which projects onto two dual modules, and $P_0 = n^{-1}\wh{\kappa}_\wedge+1$.


\section{Adjoint modules of orthogonal and symplectic groups}

Suppose now that $V$ has a nonsingular bilinear form $b$ which is either symmetric or antisymmetric, so that the group $G$ of linear transformations $V$ preserving $b$ is orthogonal or symplectic, respectively.
As a result of their geometric and physical importance, there are many techniques for the harmonic analysis of the orthogonal and symplectic groups that can be used instead of the methods employed for the unitary groups. We shall restrict ourselves to an outline of the main points without checking all the details of irreducibility, for which we refer to \cite{C,V2, LM,IK1,IK2}.

As both $b$ and $\kappa$ are nonsingular, for $u$, $v \in V$, and $X\in \LieG$, the Riesz representation theorem and universal property of tensor products tell us that there is a unique $[u\otimes v]\in \LieG$ such that the linear functional $X \mapsto b(u,X.v) = \kappa([u\otimes v],X)$.
Moreover, since $b$ is $G$-invariant we have $b(u,X.v) = -b(X.u,v)$, which is $-b(v,X.u)$ or $b(v,X.u)$, for symmetric or antisymmetric $b$, giving $[u\otimes v] = -[v\otimes u]$, or $[u\otimes v] = [v\otimes u]$, respectively. In the former (orthogonal) case this means that $[u\otimes v]$ is really a multiple of $u\wedge v$, whilst in the latter (symplectic) case it is a multiple of 
$u\otimes_S v$. Correspondingly one gets  isomorphisms $\wedge^2 V \cong \so(V)$ and $\schur^2 V \cong \symp(V)$ in the two cases, which intertwine the $G$-actions since $b$ and $\kappa$ are $G$-invariant \cite[20.1]{FH}.

The Gramian determinant provides a natural inner product  $b^{(2)}$ on $\wedge^2V$ given by 
$$
b^{(2)}(u_1\wedge u_2,v_1\wedge v_2) =
\left|\matrix{b(u_1,v_1) &b(u_1,v_2)\cr b(u_2,v_1) &b(u_2,v_2)\cr}\right|,
$$
(and for $\otimes_S^2V$ the permanent replaces the determinant). These have the same symmetries as $\kappa$ so choose the isomorphism $\wedge^2 V \cong \so(V)$ so that $\kappa = b^{(2)}$, and 
then we have
$$
\kappa(u_1\wedge u_2,v_1\wedge v_2) 
= b(u_1,b(u_2,v_2)v_1 - b(u_2,v_1)v_2),
$$
suggesting that $(v_1\wedge v_2).u_2 = b(u_2,v_2)v_1 - b(u_2,v_1)v_2$, consistent with
$b(u_1,X.u_2) = \kappa(u_1\wedge u_2,X)$ when $X = v_1\wedge v_2$, and normal conventions.
Summarising, there are  isomorphisms $\wedge^2 V \cong \so(V)$ and $\schur^2 V \cong \symp(V)$ in the two cases, which intertwine the $G$-actions since $b$ and $\kappa$ are $G$-invariant, and
thence we  have $\schur^2 \so(V) \cong \schur^2(\wedge^2V)$ and
$\schur^2 \symp(V) \cong \schur^2(\schur^2V)$.

Having shown that both cases involve submodules of $\otimes^4 V$ we shall concentrate on the orthogonal case.
The standard technique is to use harmonic tensors, which are defined as elements of the kernel of all the operators $b_{jk}: \otimes^pV \to \otimes^{p-2}V$, defined by contracting $b$ with the factors $u^{(j)}$ and $u^{(k)}$ (leaving the order of the remaining factors unaltered).
In the case of $(\wedge^2V)\otimes_S(\wedge^2V) \subset \otimes^4V$, the operators $b_{12}$ and $b_{34}$ both contract a symmetric form against an antisymmetric tensor and give zero. We can therefore concentrate on $b_{13}$, $b_{23}$, $b_{14}$, and $b_{24}$, but again the symmetries ensure that all four of these give the same result (up to sign), sending $(u\wedge v)\otimes_S(x\wedge y)$ to
$$
 b(u,x) v\otimes_S y -b(v,x) u\otimes_S y -b(u,y) v\otimes_S x + b(v,y) u\otimes_S x.
$$
Since they are essentially the same they all have the same kernel, and so harmonic tensors are those in the kernel of $b_{13}$. 

Like $\kappa$, the non-singular form $b$ has a dual form on the dual space $\wh{V}$: $\wh{b} \in \schur^2 V$, and this defines an operator $\wh{b}_{13}:\otimes^{2}V \to \otimes^4V$ by tensoring with $\wh{b}$,  but putting the new tensor factors into the first and third positions. One readily checks that $b_{13}(\wh{b}_{13}) = n = \dim(V)$, and this tells us that $p_{13} = n^{-1}\wh{b}_{13}\circ b_{13}$ is a projection. 
Clearly $\ker(p_{13}) \subseteq \ker(b_{13})$, and the construction of $p_{13}$ shows that there is equality.  The space $\schur^2(\wedge^2V)$ can be decomposed into the direct sum of the kernel and image of $p_{13}$, both invariant subspaces, where $A \in \schur^2(\wedge^2V)$ decomposes as
$$
A = (1-p_{13})A + p_{13}A.
$$
Although $b_{13}A \in \schur^2 V$ has only two factors we shall continue to label these by the original indices, 2 and 4, to avoid confusion.
It  is not harmonic but can be expresed as a sum of harmonic terms, since
\begin{eqnarray*}
p_{13}A &=& n^{-1}\wh{b}_{13}b_{13}A\\ 
&=& n^{-1}\wh{b}_{13}(1-p_{24})b_{13}A + n^{-1}\wh{b}_{13}p_{24}b_{13}A\\
&=& n^{-1}\wh{b}_{13}(1-p_{24})b_{13}A + n^{-2}\,{\wh{b}_{13}}{\wh{b}_{24}}b_{24}b_{13}A,
\end{eqnarray*}
and $(1-p_{24})b_{13}A$ is clearly harmonic, and $b_{24}b_{13}A \in \wedge^0V \cong \complex$ is a scalar, and so trivially harmonic.
Summarising, we have
$$
A = [(1-p_{13})A] + n^{-1}\wh{b}_{13}[(1-p_{24})b_{13}A] + n^{-2}\,\wh{b}_{13}\wh{b}_{24} [b_{24}b_{13}A],
$$
where all the terms in square brackets are harmonic. (This extends to a well-known  iterative procedure.) 
We already know that $b_{13}A \in \schur^2 V$, which has  dimension $\frac12n(n+1)$, and so we have effectively split $p_{13}A$ into scalars (dimension 1) and  harmonic symmetric tensors in $\schur^2V$ of dimension
$\fract12n(n+1)-1 = \fract12\left(n^2+n-2\right) = \fract12(n-1)(n+2)$,
and both components are irreducible.

Since $\schur^2(\wedge^2V)$ has dimension
\begin{eqnarray*}
\fract12\left[\frac{n(n-1)}{2}\left(\frac{n(n-1)}{2} + 1\right)\right]
&=&  \fract18n(n-1)(n^2-n +2),
\end{eqnarray*}
the first term $(1-p_{13})A$ lies in a subspace of dimension
\begin{eqnarray*}
\fract18n(n-1)(n^2-n +2) - \fract12n(n+1) 
&=& \fract18n(n-3)(n^2 + n +2),
\end{eqnarray*}
but this is not irreducible.

Since $v_1\wedge v_2\wedge v_3\wedge v_4 = v_3\wedge v_4\wedge v_1\wedge v_2$ and $(12)(34)$ are even permutations, and so $\wedge^4 V \subseteq \schur^2(\wedge^2V)$, and, since contractions with $b$ vanish, it is harmonic. Its orthogonal complement in the harmonic part of 
$\schur^2(\wedge^2V)$ has dimension 
\begin{eqnarray*}
\fract18n(n-3)\left(n^2  +n +2\right) - \fract{1}{24}n(n-1)(n-2)(n-3)
&=& \fract{1}{12}n(n+1)(n+2)(n-3).
\end{eqnarray*}
This has highest weight twice that of $\LieG = \wedge^2V$, and is the Cartan square submodule $(\wedge^2V)^{(2)} = \so(n)^{(2)}$, 
(cf \cite[7.1.2 Ex. B,C,D]{GW}).
For $n\neq 8$, the dimensions 
$$
\dim[(\wedge^2V)^{(2)}]  = \fract{1}{12}n(n+1)(n+2)(n-3), \qquad 
\dim[\wedge^4V] = \fract{1}{24}n(n-1)(n-2)(n-3),
$$
agree with those of $Y_2(\alpha)$ and $Y_2(\beta)$ given by Vogel  (when we use the orthogonal group values $\alpha = -2$, $\beta = 4$, $\gamma = n-4$, $t = n-2$ \cite[Th.1.2]{LM}, which we shall see are just minus the eigenvalues of $\wh{\kappa}$.
Heurisically, we can think of $\wedge^4V$ as $(\wedge^2V)^{(1^2)}$, since two copies of the highest weight of $\wedge^2V$ have been juxtaposed.
One can check that the decomposition for  $\so(6)$ matches that for $\su(4)$ to which it is isomorphic.
(Some insight into the geometrical significance of the decomposition can be obtained by interpreting 
$\schur^2(\wedge^2V)$ as quadratic functions on the projective space ${\mathbb P}(\wedge^2 \wh{V})$,  with $\wedge^4V$ the ideal of quadratics vanishing on Gr$_2$, \cite[15.34]{FH}.)

The above exception of $n = 8$ arises because, although $\wedge^4V$ is usually irreducible, 
it decomposes just in the case of $\so(8)$, which has the non-trivial outer automorphism group $S_3$, \cite[p33 et seq.]{A}, well known as the symmetries of its Dynkin diagram (shown on the left, with the two spin modules $\Delta^{\pm}$ ):
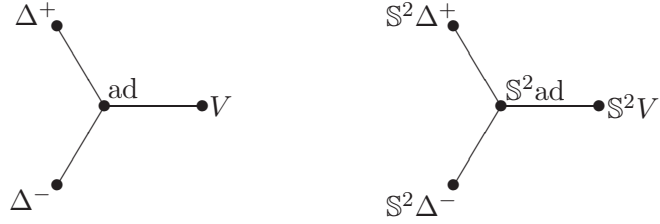
\begin{center}
{\begin{picture}(200,70)(0,-35)
\put(0,0){\line(1,0){37}}
\put(0,0){\line(-3,5){18}}
\put(0,0){\line(-3,-5){18}}
\put(0,0){\makebox(0,0){$\bullet$}}
\put(37,0){\makebox(0,0){$\bullet$}}
\put(-18,30){\makebox(0,0){$\bullet$}}
\put(-18,-30){\makebox(0,0){$\bullet$}}
\put(7,7){\makebox(0,0){$\ad$}}
\put(44,0){\makebox(0,0){$V$}}
\put(-26,35){\makebox(0,0){$\Delta^+$}}
\put(-28,-35){\makebox(0,0){$\Delta^-$}}
\put(150,0){\line(1,0){37}}
\put(150,0){\line(-3,5){18}}
\put(150,0){\line(-3,-5){18}}
\put(150,0){\makebox(0,0){$\bullet$}}
\put(187,0){\makebox(0,0){$\bullet$}}
\put(132,30){\makebox(0,0){$\bullet$}}
\put(132,-30){\makebox(0,0){$\bullet$}}
\put(163,7){\makebox(0,0){$\schur^2\ad$}}
\put(200,0){\makebox(0,0){$\schur^2V$}}
\put(120,35){\makebox(0,0){$\schur^2\Delta^+$}}
\put(120,-37){\makebox(0,0){$\schur^2\Delta^-$}}
\end{picture}}
\end{center} 
On the right is the corresponding diagram for the symmetric squares. It is known, \cite[Th. 4.6(i)]{A} that $\schur^2\Delta^\pm = (\wedge^4V)^\pm \oplus \complex$, where $(\wedge^4V)^\pm$ denote the $\pm1$-eigenspaces of the Hodge star operator, whose direct sum is $\wedge^4V$. We also have $\schur^2V = \schur_b^2V \oplus \complex$, and, on setting $n=8$, we find that $\dim(\schur_b^2V)=35$.
By the outer automorphism symmetries these three modules must have the same dimension, and 
$\wedge^4V$, which,  by our earlier formulae, has dimension 70, splits into two 35-dimensional eigenspaces of the Hodge star operator.

We could find the values of the Casimir operators on the irreducible submodules by linking the Casimir operator to the Laplacian on the underlying spaces of symmetric functions, but it is simpler to use the known highest weights under a maximal torus $T$.
Taking an orthonormal basis $\{\epsilon_j: j=1,\ldots,\lfloor\frac12n\rfloor\}$ of $\LieT$, the highest root of $\LieG$ is $\lambda = \epsilon_1+\epsilon_2$, and the sum of the positive roots is given by $2\delta = \sum_j (n-2j)\epsilon_j$, \cite[Planches II et IV]{B}
This gives $\wh{\kappa}(\lambda,\lambda) = 2$, and $\wh{\kappa}(\lambda,2\delta) = (n-2) +(n-4) = 2(n-3)$, so that the quadratic Casimir takes the value $2(n-2)$.
The Cartan square with highest weight $2\lambda$ now has Casimir with value $4\times 2 + 4(n-3) = 4(n-1)$.
The highest weight for $\wedge^4V$ is $\lambda = \epsilon_1+\epsilon_2+\epsilon_3+\epsilon_4$, giving $\wh{\kappa}(\lambda, \lambda) = 4$ and $\wh{\kappa}(\lambda,2\delta) = 4(n-5)$. The Casimir therefore takes the value $4(n-4)$.
Finally the highest weight in $\schur^2V$ is $2\epsilon_1$, which gives 
$\wh{\kappa}(\lambda, \lambda) = 4$ and $\wh{\kappa}(\lambda,2\delta) = 2(n-2)$, so that $C$ takes the value $2n$.
  
We can now calculate the eigenvalues of the polarised Casimir $\wh{\kappa}_\schur$ on $\schur^2\LieG$ to be 
 \begin{eqnarray*}
 {\LieG^{(2)}}:  2(n-1) - 2(n-2) = 2,  &\qquad& 
 {\wedge^4V}: 2(n-4) - 2(n-2) = -4,\\
{\schur_b^2V}: n- 2(n-2) = 4-n &\qquad&
 {\complex}: 0 - 2(n-2) = 4-2n.
\end{eqnarray*}
These are again just minus the Vogel parameters. This shows the characteristic equation of $\wh{\kappa}_\schur$ on $\schur^2\LieG$ to be
$$
(\wh{\kappa}_\schur -2)(\wh{\kappa}_\schur +4)(\wh{\kappa}_\schur +n-4)(\wh{\kappa}_\schur +2n-4) = 0,
$$
from which the projections into eigenspaces can be found as for unitary groups:
\begin{eqnarray*}
P_{2} = \frac{(\wh{\kappa}_\schur +4)(\wh{\kappa}_\schur +n-4)(\wh{\kappa}_\schur +2n-4)}{12(n-2)(n-1)}, &\quad&
P_{-4} = \frac{(\wh{\kappa}_\schur -2)(\wh{\kappa}_\schur +n-4)(\wh{\kappa}_\schur +2n-4)}{-12(n-8)(n-4)}, \\
P_{4-n} = \frac{(\wh{\kappa}_\schur -2)(\wh{\kappa}_\schur +4)(\wh{\kappa}_\schur +2n-4)}{n(n-2)(n-8)} &\quad&
P_{4-2n} = \frac{(\wh{\kappa}_\schur -2)(\wh{\kappa}_\schur +4)(\wh{\kappa}_\schur +n-4)}{-4n(n-1)(n-4)}.
\end{eqnarray*}

We can similarly investigate $\wedge^2 \LieG$. We already know that there is a copy of the adjoint module, but there is only one other irreducible submodule (since $b$ merges the two dual submodules in the unitary case).
Since $\wedge^2 \LieG$ has dimension
$$
\fract12\left[\frac{n(n-1)}{2}\left(\frac{n(n-1)}{2} - 1\right)\right]
=  \fract18n(n-1)(n^2-n -2),
$$
the other irreducible must have dimension
\begin{eqnarray*}
\fract18n(n-1)(n^2-n -2)- \fract12n(n-1) 
&=& \fract18n(n-1)(n-3)(n+2).
\end{eqnarray*}
Higher tensor powers can be treated similarly by extending the harmonic decomposition.

In the case of symplectic groups the symmetric and exterior products such as $\schur^2 V$ and $\wedge^2 V$  are interchanged, and  the dimension $\frac12n(n+1)$ of the former is replaced by  $\frac12n(n-1)$ for the latter. Formally this is the same as changing the sign of $n$ or replacing $V$ by the virtual $-V$, and that also works when we replace $\wedge^4V$ by $\schur^4V$.
We therefore find the dimensions of the irreducible components, to be 1, $\frac18(n^2-n-2)$, $\frac1{24}n(n+1)(n+2)(n+3)$, for $\schur^4V$, and $\fract{1}{12}n(n-1)(n-2)(n+3)$ for $\LieG^{(2)}$. 

The above examples of the classical unitary, orthogonal and symplectic groups can be extended to the exceptional groups $G_2$, $F_4$, $E_6$, $E_7$, and $E_8$, but, even exploiting some of the shortcuts in \cite{A} each needs individual treatment.
As noted in the Introduction, $\wedge^2\LieG$ contains a copy of $\LieG$, and 
$\schur^2\LieG$ contains the trivial module and the Cartan square $\LieG^{(2)}$.
It turns out that in each case there is just one more irreducible.  Since the dimensions of 
$\schur^2\LieG$ and $\wedge^2\LieG$ are easily calculated we only need to calculate the dimension of the Cartan square to know the dimensions of all the irreducible components.
The results are as follows

\begin{center}
\begin{tabular}{|l l | l l|}\hline
$\schur^2\LieG$  &$1 \oplus \LieG^{(2)} \oplus *$ &$\wedge^2\LieG$ &$\LieG \oplus *$\\ \hline\hline
$\schur^2\LieG_2$  &$1 \oplus 77 \oplus 27 $ &$\wedge^2\LieG_2$ &$14 \oplus 77^*$\\
$\schur^2\LieF_4$  &$1 \oplus 1053 \oplus  324$ &$\wedge^2\LieF_4$ &$52 \oplus 1274$\\
$\schur^2\LieE_6$  &$1 \oplus 2430 \oplus 650$ &$\wedge^2\LieE_6$ &$78 \oplus 2925$\\
$\schur^2\LieE_7$  &$1 \oplus 7371 \oplus 1539$ &$\wedge^2\LieE_7$ &$133 \oplus 8645$\\
$\schur^2\LieE_8$  &$1 \oplus 27000 \oplus 3875$ &$\wedge^2\LieE_8$ &$248 \oplus 30380$\\ \hline
\end{tabular}
\end{center}

\medskip 
For $\LieG_2$ one can quite easily find the dimension of the Cartan square using Freudenthal's construction, \cite[\S22.2, \S22.4]{FH}.
This realises $\LieG_2$ as the direct sum of $\LieG_0 = \su(3)$ with its natural module $V=\complex^3$ and the natural dual $\wh{V}$ , equipped with the Lie bracket on $\su(3)$, the adjoint action of $\su(3)$ on $V$ and $\wh{V}$ being the natural actions, the Lie bracket on $V$ being given by 
$V\wedge V \cong \wh{V}$, and the dual on $\wh{V}$, and the Lie bracket of $V$ with $\wh{V}$, coming from the antisymmetrised identification of $V\otimes \wh{V}$ with the trace free part of $\hom(V)$.

The  Cartan subalgebra $\LieT$ of $\su(3)$ can also be used for $\LieG_2$. In addition to the six root vectors in $\su(3)$ there are three more (weight vectors) in $V$ and also  in $\wh{V}$. Identifying the dual $\wh{\LieT}$ with elements of $\real^3$ whose components add up to 0, we can take the roots in $V$ to be   $\alpha_1 = (-1,0,1)$ and its cyclic permutations $\alpha_2 = (1,-1,0)$ and $\alpha_3 = (0,1,-1)$, and those in $\wh{V}$ are their negatives. We take the positive roots to be $\alpha_1$, $\alpha_2$ in $V$, and 
$-\alpha_3 \in \wh{V}$, together with the positive roots $\alpha_1 - \alpha_2 = (-2,1,1)$, 
$\alpha_2 - \alpha_3 = (1,-2,1)$ and $\alpha_1 -\alpha_3 = (-1,-1,2)$ in $\su(3)$, with the last being the highest root $\lambda$ in both $\su(3)$ and $\LieG_2$.
The sum of the positive roots is $2\delta = 2(-1,-2,3)$.

Weyl's dimension formula gives the dimension as the product over positive roots $\alpha$ of 
$1+(\lambda,\alpha)/(\delta,\alpha)$ and, in our case, we need not include $\alpha_2$ which is orthogonal to $\lambda$.  The inner product $(\lambda,\lambda) = 6$, and all the other 
$(\lambda, \alpha) $ are 3. One can similarly find the inner products $(\delta,\alpha)$, and, as a check, the dimension comes out as
$$
\left(1 +\frac33\right)\left(1 +\frac34\right)\left(1 +\frac35\right)\left(1 +\frac36\right)\left(1 +\frac69\right) = 2.\frac74.\frac85.\frac32.\frac53 = 14 
$$
We find the dimension of the Cartan square or cube, by changing $\lambda$ to  $2\lambda$ or $3\lambda$, respectively:
\begin{eqnarray*}
\LieG_2^{(2)}&:& \left(1 +\frac63\right)\left(1 +\frac32\right)\left(1 +\frac65\right)\left(1 +\frac66\right)\left(1 +\frac{12}9\right) = 3.\frac52.\frac{11}5.2.\frac73 = 77.,\\
\LieG_2^{(3)}&:& \left(1 +\frac93\right)\left(1 +\frac94\right)\left(1 +\frac95\right)\left(1 +\frac96\right)\left(1 +\frac{18}9\right) = 4.\frac{13}4.\frac{14}5.\frac52.3 =  273.
\end{eqnarray*}

A similar procedure can be used for $\LieE_8$  starting with $\su(9)$ and its two natural modules $V = \complex^9$, and its dual, but we model $\LieE_8$ on $\su(9)\oplus \wedge^3V\oplus \wedge^3\wh{V}$ which is $80 +84 +84 = 248$ dimensional. The obvious exterior product map 
$(\wedge^3V\otimes \wedge^3V \to \wedge^6V \cong \wedge^3\wh{V}$ helps us to proceed as before.
Unfortunately, we must work with 120 positive roots, and even though many of those are orthogonal to the highest root, we are still left with almost 60 terms to compute.
Another commonly used approach is through the Witt construction, which is similar to that of Freudenthal, but with vanishing Lie brackets of two elements in $M$ or in $\wh{M}$. One can then use the 120-dimensional $\LieG_0 = \spin(16)$ with its 128-dimensional even spinor module $\Delta^+$, but this still requires similar intense calculation to find the dimension of the Cartan square.

\appendix


\section{Weyl's dimension formula}

In this Appendix we want to show that for special unitary groups the various submodules which we found are irreducible, by showing that they are the highest weight modules with their particular highest weights. Although the calculations can be done by computer, it is instructive to see the 
general expressions which apply to all the highest weights relevant to our examples.

Weyl's dimension formula gives the dimension of the irreducible SU($n$)-module with highest weight $\lambda = (\lambda_1,\ldots,\lambda_n)$ as, \cite[(6.3(1)]{FH}, \cite[(7.1.17)]{GW}:
$$
\prod_{i<j}\left(\frac{\lambda_i - \lambda_j + j - i}{j-i}\right).
$$
Clearly, we need only consider the terms where $\lambda_i \neq \lambda_j$, and we can take the product of the cases where $\lambda_i - \lambda_j  = 1,2,\ldots$, and we now investigate regions where these differences are constant.

The simplest nontrivial case $\Phi_1(a,b,c,d)$ occurs where $\lambda_i - \lambda_j = 1$ for $a+1\leq i\leq b$ and $c+1\leq j\leq d$ with $a < b \leq c < d$, where we find that
\begin{eqnarray*}
\Phi_1(a,b;c,d) 
&=& \frac{(d-a)!}{(d-b)!}\frac{(c-b)!}{(c-a)!}.
\end{eqnarray*}
This \lq\lq cross ratio\rq\rq\ expression simplifies in various special cases which will be useful later:
$$
\Phi_1(a,a+1;c,d) = \frac{d-a}{c-a}, \qquad 
\Phi_1(a,b;d-1,d)= \frac{d-a}{d-b}, \qquad
 \Phi_1(a,b;b,d) = {{d-a}\choose{b-a}}.
$$

Where the $\lambda$ components differ by 2, we need
\begin{eqnarray*}
\Phi_2(a,b;c,d) :=
\prod_{j=c+1}^{d}\prod_{i=a+1}^{b} \frac{j-i+2}{j-i} 
= \Phi_1(a,b;c+1,d+1)\Phi_1(a,b;c,d),
\end{eqnarray*}
and, more generally, for regions where the $\lambda$ components differ by $m$ we similarly have
$$
\Phi_m(a,b;c,d) = \prod_{r=0}^{m-1} \Phi_1(a,b;c+r,d+r),
$$
and, in particular, 
$$
\Phi_m(0,1;c,d) = \prod_{r=0}^{m-1}\frac{d+r}{c+r}
= \frac{(d+m-1)!(c-1)!}{(d-1)!(c+m-1)!}
= \Phi_1(1-m,1;c,d),
$$
and, similarly, $\Phi_m(a,b;c,c+1) = \Phi_1(a,b;c,c+m)$.
When both constraints hold  this reduces further to 
$\Phi_m(0,1;d-1,d) = (d+m-1)/(d-1) = (c+m)/c$.

We start by considering the case when the highest weight has the form $(1^{b-a},0^{c-b},(-1)^{d-c})$, where the superscripts denote the number of repetitions.  This produces differences of 1 between the middle and outer regions, and 2 between the two outer regions, so that the dimension is given by 
$$
\Psi_1(a,b;c,d) = \Phi_1(a,b;b,c)\Phi_1(b,c;c,d)\Phi_1(a,b;c,d)\Phi_1(a,b;c+1,d+1).
$$
The first three terms give a trinomial coefficient $(d-a)!/[(d-c)!(c-b)!(b-a)!]$, 
and $\Psi_1$ is
\begin{eqnarray*}
\frac{(d-a)!}{(d-c)!(c-b)!(b-a)!}\frac{(d-a+1)!}{(d-b+1)!}\frac{(c-b+1)!}{(c-a+1)!}
&=& \frac{c-b+1}{d-a+1}{{d-a+1}\choose{d-c}}{{d-a+1}\choose{b-a}},
\end{eqnarray*}
a formula with some interesting applications.

The highest weight in $\LieG^{(1^k)} \subset \hom(\wedge^kV)$ when $\dim(V) = n$ is $(1^k,0^{n-2k},(-1)^k)$.
In this case  
$$
\dim(\LieG^{1^k)} = \frac{(n-2k+1)}{n+1}{{n+1}\choose{k}}^2.
$$

The highest weight in $\LieG^{(k)} \subset \hom(\schur^kV)$ is $(k,0^{n-2},-k)$ and the dimension is given by 
\begin{eqnarray*}
\dim(\LieG^{(k)} &=& (n+2k-1)/(n-1){{n+k-2}\choose{k}}^2
\end{eqnarray*}

The irreducible components of $\otimes^2(\wedge^kV)$ have highest weights 
$(2,\ldots,2,1,\ldots,1,0\ldots,0)$ with changes after the positions $p=k-r$ and $q=k+r$, for 
$0 \leq  r < \min(k,n-k)$, 
and our formula gives the dimension as
$$
\frac{(2r+1)}{n+1}{{n+1}\choose{k-r}}{{n+1}\choose{n-k-r}}.
$$
When $r$ is odd this provides an irreducible component of $\wedge^2(\wedge^kV)$, and when $r$ is even an irreducible component of  $\schur^2(\wedge^kV)$, in particular, when $r = 0$ we get the Cartan square  with dimension 
$$
\dim((\wedge^kV)^{(2)}) = \frac{1}{n+1}{{n+1}\choose{k}}{{n+1}\choose{k+1}}.
$$
When $k=2$ this reduces to $\frac16{{n^2}\choose {2}}$, as mentioned in Section 4.

Other highest weights such as $(1^k,0,\ldots,-k)$ and $(k,0,\ldots,(-1)^k)$ which appear in the antisymmetric tensor product can be dealt with in a similar way. 
(There is a highest weight vector for $\LieG^{(1^k,-k)} \subset \hom(\schur^k V,\wedge^k V)$ with weight $(1^k,0,\ldots,-k)$. 
Similarly in $\hom(\wedge^m V,\schur^m V)$ one finds a highest weight $(k,0,\ldots,(-1)^k)$
associated with the dual  irreducible submodule $\LieG^{(k,(-1)^k)}$.)
This module, like its dual,  has dimension
\begin{eqnarray*}
\dim(\LieG^{(1^k,-k)}) &=& {{n-1}\choose{k}}{{n+k}\choose{k}} = \frac{\prod_{r=1}^k (n^2-r^2)}{(k!)^2}.
\end{eqnarray*}

The same functions appear in Weyl's dimension formula for orthogonal groups although it is rather more complicated:
$$
\prod_{i < j} \frac{(\lambda_i +\delta_i)^2 - (\lambda_j +\delta_j)^2}{\delta_i^2 - \delta_j)^2}
= \prod_{i < j} \frac{(\lambda_i - \lambda_j)+(\delta_i - \delta_j)}{\delta_i - \delta_j}
\prod_{i < j} \frac{\lambda_i  + \lambda_j +\delta_i+\delta_j}{\delta_i + \delta_j}
$$
where the indices run from 1 to $\ell$, the rank of the orthogonal group, and $\delta_j = \ell - j$ for all 
$j$.
We shall be interested in the case  when $\lambda_i = m$ for $i\leq b$, and vanishes for $i > b$, and, for definiteness, concentrate on the even case where $n = 2\ell$.
Conveniently, we now have $\lambda_i \pm \lambda_j = \lambda_i$, when $j >b$, and, as in the unitary case, we only need the factors where $\lambda_i \neq 0$.
This leaves us with 
$$
\prod_{i =1}^b  \prod_{j = b+1}^\ell \frac{m +j-i}{j - i}
\prod_{i =1}^b  \prod_{j = b+1}^\ell  \frac{m + n - (j+i)}{ n - (j+i)}
\prod_{1\leq i <j\leq b}\frac{(2m + n - i -j)}{(n - i - j)} .
$$
The leading double product is $\Phi_m(0,b;b,\ell)$, and, setting $k = n-j$,   the second  becomes 
\begin{eqnarray*}
\Phi_m(0,b;\ell-1,n-b-1) 
&=& \Phi_m(0,b;\ell,n-b-1)\prod_{i=1}^b\frac{(2m + n - 2i)}{(n -2 i)}.
\end{eqnarray*}
The extra term at the end combines with that at the end of the earlier expression, and using 
 $\Phi_m(0,b;b,\ell)\Phi_m(0,b;\ell,n-b-1) = \Phi_m(0,b;b,n-b-1)$,  the product 
reduces  to 
$$
\Phi_m(0,b;b,n-b-1)\prod_{1\leq i \leq j\leq b}\frac{(2m + n - i -j)}{(n - i - j)}
$$
(For odd-dimensional orthogonal groups an extra product term precisely compensates for the disappearance of the overlap between the two double products and matches that for even groups.)

This formula is easily applied to the cases of $\wedge^2V \cong \so(n)$, $\so(n)^{(1^2)}$, and the Cartan square $\so(n)^{(2)}$ as tabulated below:

\begin{center}
\begin{tabular}{|l|l|l|c|}\hline
SO($n$)-module &$m$ &$b$ &dimension\\ \hline
 $\so(n)$ &1 &2 &${{n}\choose{2}}$ \\
 $\so(n)^{(1^2)}$ &1 &4 &${{n}\choose{4}}$\\
 $\so(n)^{(2)}$  &2 &2 &$\fract12(n-3){{n+2}\choose{3}}$\\
 \hline
\end{tabular}
\end{center}
\smallskip
in  agreement with Section 4.


\section{Highest weight vectors}

Highest weight vectors are usually defined in terms of a chosen maximal torus (or, more accurately, borel subalgebra), but all choices are known to be equivalent for their main applications. 
Another definition of highest weight vectors is independent of any choice  \cite{KCH82,WL82}, and provides a criterion is simple enough to make sense for other groups and Hopf algebras. Applied to the Heisenberg group it picks out coherent states, for Clifford algebras the pure spinors, and for loop groups the tau functions for integrable systems, \cite{KCH00}.

\medskip{\it
Let $V$ be an irreducible module for a compact simple Lie group $G$ (and its Lie algebra $\LieG$).
A vector $v$ in $V$ is a highest weight vector for some choice of maximal torus/borel subalgebra if and only if $v^{\otimes2} = v\otimes v$ generates an irreducible submodule of the tensor product.
}

\medskip
To prove this result, we need to compare this definition with the usual one. If $v$ is a highest weight vector for some choice of $T$ then $v^{\otimes2} = v\otimes v$ is a highest weight for the Cartan square $V^{(2)}$, which is an irreducible submodule of $V\otimes V$, proving the necessity of that condition. 
We note that $V$ is necessarily finite dimensional, with an invariant inner product, and both these properties extend to $V\otimes V$. To prove sufficiency we suppose that $u \in V$ is a highest weight vector for some choice of maximal torus $T$ (Cartan subalgebra $\LieT$).
Since $u$ generates the irreducible $V$ we know that $\ip{v}{gu}$ cannot vanish identically for all $g\in G$, so neither can $\ip{(g^{-1}\!.v)^{\otimes2}}{u^{\otimes2}} = \ip{v}{g.u}^2$, and so the submodule generated by $v^{\otimes2}$ contains $u^{\otimes2}$, the highest weight vector for the Cartan square $V^{(2)}$ (the irreducible with highest weight twice the highest weight for $V$), and the module generated by $v^{\otimes2}$ must contain $V^{(2)}$. Consequently if $v^{\otimes2}$ generates a single irreducible submodule of $V\otimes V$, that irreducible is $V^{(2)}$. (This argument extends to higher tensor powers, and  if $v^{\otimes m}$ generates an irreducible submodule then that irreducible must be the $m$-th  Cartan power, $V^{(m)}$.)

More detailed information about $v$ is provided by the polarised Casimir operator $\wh{\kappa}$, which intertwines the action of $G$ on the irreducible $V^{(2)}$ and acts as multiplication by 
$\wh{\kappa}(\lambda,\lambda)$, so that $\wh{\kappa}v^{\otimes2} = \wh{\kappa}(\lambda,\lambda)v^{\otimes2}$.
Recalling the first paragraph of Section 3, 
we introduce $H_v = \ip{v}{\wh{\kappa}_1.v}\wh{\kappa}_2 \in \LieG$, so, for normalised $v$, the eigenvalue property for $\wh{\kappa}$ shows that $v$ is an eigenvector for $H_v$: 
$$
H_v v = \ip{v}{\wh{\kappa}_1.v}\wh{\kappa}_2.v =  \wh{\kappa}(\lambda,\lambda)v.
$$
We also note that for any $Y\in \LieG$ we can use the expansion formula 
$Y = \wh{\kappa}_1\kappa(\wh{\kappa}_2,Y)$ to obtain 
$$
\kappa(H_v,Y) = \ip{v}{\wh{\kappa}_1v}\kappa(\wh{\kappa}_2,Y) = \ip{v}{Yv}.
$$
In particular, 
taking $Y= H_v$ and normalising $v$,  we have
$$
\kappa(H_v,H_v) = \ip{v}{H_v.v} = \wh{\kappa}(\lambda,\lambda).
$$ 

Since $H_v \in \LieG$,  it is conjugate to an element of $\LieT$, that is $gH_v g^{-1} \in \LieT$, and, moreover, the invariance of $\wh{\kappa}$ shows that  $gH_v g^{-1}= H_{gv}$.
We shall write $w = gv$.
We can expand $H_w$ in terms of an orthonormal basis $\{H_j;j = 1,2,\ldots,r\}$ of the Cartan subalgebra $\LieT$ as $H_w  = \sum_j \alpha_jH_j$, where, comparing norms, we have 
$$
\wh{\kappa}(\lambda,\lambda) = \wh{\kappa}(\alpha,\alpha).
$$
One can also expand $w$ in terms of weight vectors $u_\mu$, and, since $w$ is an eigenvector of 
$H_w$, the expansion must only involve weights $\mu$ with 
$\wh{\kappa}(\lambda,\lambda) = \mu(H_w) = \wh{\kappa}(\mu,\alpha)$. Since $\wh{\kappa}(\mu,\mu) \leq \wh{\kappa}(\lambda,\lambda)$ we have
$$
\wh{\kappa}(\lambda,\lambda)^2 = \wh{\kappa}(\mu,\alpha)^2 \leq 
\wh{\kappa}(\mu,\mu)\wh{\kappa}(\alpha,\alpha) 
\leq \wh{\kappa}(\lambda,\lambda)^2,
$$
which is only possible with saturated inequalities:  
$\wh{\kappa}(\mu,\mu) = \wh{\kappa}(\lambda,\lambda)$ and $\wh{\kappa}(\mu,\alpha)^2 = \wh{\kappa}(\mu,\mu)\wh{\kappa}(\alpha,\alpha)$.
The former happens only when $\mu$ is on the Weyl group orbit of $\lambda$ :
That means that for some $n \in N(T)$ we have
 $\mu(H) = (\ad^*(n^{-1})\lambda)(H) = \lambda(nHn^{-1})$ for all $H \in \LieT$.
 But, for any $H \in \LieT$ we have 
$$
H.(n.u) = n(n^{-1}Hn).u = n\mu(n^{-1}Hn).u  = \lambda(H)(n.u).
$$
so that $u = n.w = ng.v$ is the highest weight vector for $T$, and $v$ the highest weight vector for 
$(ng)T(ng)^{-1}$, as asserted in the theorem.

\medskip
\noindent{\it Remark.}
Writing $H_{nw}  = \sum_j \beta_jH_j$, our previous arguments show that $\wh{\kappa}(\beta,\beta) = 
\wh{\kappa}(\lambda,\lambda)$ 
and saturation of the Cauchy--Schwarz inequality gives $\beta = \lambda$ and $H_u  = \sum_j \lambda_jH_j$.
Using $H_u$ in an earlier identity, we have, for any $Y \in \LieT$,
$$
\kappa(H_u,Y)  = \ip{u}{Y.u} = \lambda(Y),
$$
so that $H_u$ is the vector representing $\lambda$ in the Riesz representation theorem.

In the case of the adjoint module, where $Y:  v \mapsto [Y,v]$, $\lambda=\alpha$ is the highest root, and the Killing form and inner product are connected by defining an antilinear involution $v \mapsto v^*$ such that $\ip{v}{u} = \kappa(v^*,u)$, we have $H_v  = \sum_j \kappa(v^*,[\wh{\kappa}_1,v])\wh{\kappa}_2   = \sum_j \kappa([v,v^*],\wh{\kappa}_1)\wh{\kappa}_2 = [v,v^*]$.
It also follows from the equation $[H_v,v] = \wh{\kappa}(\alpha,\alpha)v$ that $\kappa(v,v) = 0$, so that $v$ is nilpotent. 
Now set $e_\alpha = cv$, where $c \in \real$, $f_\alpha = c v^*$, and  $h_\alpha :=[e_\alpha,f_\alpha] = c^2 H_v$. Then the equation $[H_v,v] = \wh{\kappa}(\alpha,\alpha)v$ becomes $[h_\alpha,e_\alpha] = c^2\wh{\kappa}(\alpha,\alpha)e_\alpha$ (and thence $[h_\alpha,f_\alpha] = -c^2\wh{\kappa}(\alpha,\alpha)f_\alpha$), so that one can recover the usual relations by taking $c = \sqrt{2/\wh{\kappa}(\alpha,\alpha)}$, which also gives the customary relation $H_u = \frac12\wh{\kappa}(\alpha,\alpha)\,h_\alpha$ between the Riesz representative of $\alpha$ and the coroot $h_\alpha$.

\bigskip
Results are more elusive when $v^{\otimes 2}$ generates the direct sum of several irreducible submodules of the symmetric tensor product module $\schur^2V$ (one of which must be the Cartan square submodule $V^{(2)}$). Unfortunately,  most vectors are not weight vectors for any maximal torus.  (There are only a finite number of weight vectors for a given torus $T$, and the orbit of each has dimension $\dim(G/T)$ and is of measure zero in any module with larger dimension.)
We therefore consider only the case when $v$ is a weight vector with weight $\theta$ for some maximal torus $T$. (Otherwise, the decomposition is determined by the lowest weight component of $v$.)

Nonetheless, some information is provided by the Frobenius Reciprocity Theorem, which tells us that 
for a compact group $G$, and closed  subgroup $H \rmap{i} G$, the restriction functor $i^*$ taking 
$G$-modules to $H$-modules, and the functor $i_*$ inducing $H$-modules to $G$-modules are adjoints.  More precisely, for $\hom_G$ and $\hom_H$ maps intertwining $G$-modules and $\hom_H$-modules, respectively, we have, for any $H$-module $\Theta$, and $G$-module $\Gamma$,
$$
\hom_H(i^*\Gamma,\Theta) \cong \hom_G(\Gamma,i_*\Theta). 
$$
The isomorphism takes an $H$ intertwining operator $\Lambda$ on the left to the $G$-intertwining operator $\Lambda^*$ on the right, defined by $(\Lambda^*m)(g) = \Lambda(g.m)$ for all $m \in \Theta$. Since $\Lambda$ intertwines the $H$-actions we obtain the induced function condition 
$(\Lambda^*m)(hg) = h.(\Lambda(g.m)) = h.(\Lambda^*m)(g)$, and the $G$-intertwining condition $(\Lambda^*(x.m))(g) = \Lambda(g.(x.m)) = (\Lambda^*m)(gx)$.
The inverse map recovers $\Lambda m$  by evaluating $\Lambda^*m$ at the identity.
(Further details can be found in \cite{B65}, including a geometrical interpretation of Weyl's character formula and its link to Dirac operators.)

We can apply this to the case where $H=T$ the maximal torus in  a compact simple Lie group $G$, 
$\Theta = v^{\otimes 2}$, and  $\Lambda$ takes $m$ in the $G$-submodule $M \subseteq \schur^2V$ generated by $v^{\otimes 2}$ to its projection onto $v^{\otimes2}$ ($\Lambda m = v^{\otimes2}\ip{v^{\otimes2}}{m}$, when $v$ is normalised).
The Reciprocity Theorem tells us that the induced $G$-module $i^*\Theta$ contains an irreducible 
$\Gamma$ as often as $i^*\Gamma$ contains $\Theta$, which imposes the constraint that only irreducibles large enough to contain vectors of weight $\theta^2$ can occur in the submodule generated by $v^{\otimes 2}$.

The decomposition of $\schur^2V$ imposes an upper limit on the size of the irreducibles, since 
it can contain no irreducibles with weights higher than the square of  the highest weight $\lambda$ for $V$. A complementary upper constraint comes since only direct sumands of $\schur^2V$ can occur, which rules out anything with a highest weight larger than $\lambda^2$.
Moreover any irreducible satisfying both constraints will be present.
For example, when $v$ is a vector of weight $m$ in the irreducible $D^\ell$ of $G = $SO(3) with highest weight $\ell$ we find that the module generated by $v^{\otimes 2}$ is $D^{2\ell} \oplus D^{2(\ell-1)}\oplus \ldots \oplus D^{2|m|}$. (When $m = 0$ one has the module generated by $v^{\otimes 2}$ is $D^{2\ell} \oplus D^{2(\ell-1)}\oplus \ldots \oplus D^0 = \schur^2 D^\ell$.)





\end{document}